\documentclass[12pt,reqno]{amsart}
\usepackage{amsmath , amssymb }
\usepackage{graphics}

\usepackage[mathscr]{eucal}

\newdimen\unit\newdimen\psep\newcount\nd\newcount\ndx\newbox\dotb\newbox\ptbox
\newdimen\dx\newdimen\dy\newdimen\dxx\newdimen\dyy\newdimen\hgt
\newdimen\xoff\newdimen\yoff
\newcommand\clap[1]{\hbox to 0pt{\hss{#1}\hss}}
\newcommand\vdisk[1]{{\font\dotf=cmr10 scaled #1\dotf.}}
\newcommand\varline[2]{\setbox\dotb\hbox{\vdisk{#1}}\xoff=-.5\wd\dotb
\wd\dotb=0pt\yoff=-.5\ht\dotb\psep=#2\ht\dotb}
\newcommand\varpt[1]{\setbox\ptbox\clap{\vdisk{#1}}\setbox\ptbox
\hbox{\raise-.5\ht\ptbox\box\ptbox}}
\newcommand\cpt{\copy\ptbox}
\newcommand\point[3]{\rlap{\kern#1\unit\raise#2\unit\hbox{#3}}}
\newcommand\setnd[4]{\dx=#3\unit\advance\dx-#1\unit\divide\dx by\psep
\dy=#4\unit\advance\dy-#2\unit\divide\dy by\psep \multiply\dx
by\dx\multiply\dy by\dy\advance\dx\dy\nd=1\advance\dx-1sp
\loop\ifnum\dx>0\advance\dx-\nd sp\advance\nd1\advance\dx-\nd
sp\repeat}
\newcommand\dl[4]{{\setnd{#1}{#2}{#3}{#4}\dline{#1}{#2}{#3}{#4}\nd}}
\newcommand\dline[5]{{\nd=#5\hgt=#2\unit\dx=#3\unit\advance\dx-#1\unit
\divide\dx by\nd\dy=#4\unit\advance\dy-#2\unit\divide\dy by\nd
\advance\hgt\yoff\rlap{\kern#1\unit\kern\xoff\loop\ifnum\nd>1\advance\nd-1
\advance\hgt\dy\kern\dx\raise\hgt\copy\dotb\repeat}}}

\newcommand\qellip[4]{{\setnd{0}{0}{#3}{#4}\dx=\unit\dy=0pt\raise\yoff\rlap{%
\kern#1\unit\kern\xoff\raise#2\unit\hbox{\loop\ifnum\dx>0\rlap{\kern#3\dx
\raise#4\dy\copy\dotb}\hgt=\dx\divide\hgt
by\nd\advance\dy\hgt\hgt=\dy \divide\hgt
by\nd\advance\dx-\hgt\repeat\rlap{\raise#4\dy\copy\dotb}}}}}
\newcommand\bez[6]{{\setnd{#1}{#2}{#3}{#4}\ndx=\nd\setnd{#3}{#4}{#5}{#6}
\ifnum\ndx>\nd\nd=\ndx\fi\dx=#3\unit\advance\dx-#1\unit\dy=#4\unit
\advance\dy-#2\unit\dxx=#5\unit\advance\dxx-#1\unit\dyy=#6\unit\advance
\dyy-#2\unit\advance\dxx-2\dx\advance\dyy-2\dy\divide\dxx
by\nd\divide\dyy
by\nd\advance\dx.25\dxx\advance\dy.25\dyy\divide\dx
by\nd\divide\dy by\nd \multiply\nd
by2\dx=100\dx\dy=100\dy\dxx=100\dxx\dyy=100\dyy\divide\dxx by\nd
\divide\dyy
by\nd\hgt=#2\unit\raise\yoff\rlap{\kern#1\unit\kern\xoff
\raise\hgt\copy\dotb\loop\ifnum\nd>0\advance\nd-1\advance\hgt0.01\dy
\kern0.01\dx\raise\hgt\copy\dotb\advance\dx\dxx\advance\dy\dyy\repeat}}}
\newcommand\ptu[3]{\point{#1}{#2}{\cpt\raise1ex\clap{$\scriptstyle{#3}$}}}
\newcommand\ptd[3]{\point{#1}{#2}{\cpt\raise-1.8ex\clap{$\scriptstyle{#3}$}}}
\newcommand\ptr[3]{\point{#1}{#2}{\cpt\raise-.4ex\rlap{$\ \scriptstyle{#3}$}}}
\newcommand\ptl[3]{\point{#1}{#2}{\cpt\raise-.4ex\llap{$\scriptstyle{#3}\ $}}}
\newcommand\ptlu[3]{\point{#1}{#2}{\raise.8ex\clap{$\scriptstyle{#3}$}}}
\newcommand\ptld[3]{\point{#1}{#2}{\raise-1.6ex\clap{$\scriptstyle{#3}$}}}
\newcommand\ptlr[3]{\point{#1}{#2}{\raise-.4ex\rlap{$\,\scriptstyle{#3}$}}}
\newcommand\ptll[3]{\point{#1}{#2}{\raise-.4ex\llap{$\scriptstyle{#3}\,$}}}
\newcommand\pt[2]{\point{#1}{#2}{\cpt}}

\newcommand\medline{\varline{800}{.5}}
\newcommand\thnline{\varline{400}{.6}}

\varpt{2500}\thnline\unit=2em

\newtheorem{thm}{Theorem}

\newtheorem{conj}{Conjecture}

\newtheorem{prob}{Problem}
\newtheorem{lemma}[thm]{Lemma}
\newtheorem{cor}[thm]{Corollary}

\newtheorem{obs}[thm]{Observation}
\theoremstyle{definition}\newtheorem{rmk}{Remark}
\theoremstyle{definition}\newtheorem*{defn}{Definition}
\newcommand{\ds}{\displaystyle}
\newcommand{\ul}{\underline}

\def\N{\mathbb{N}}

\def\Z{\mathbb{Z}}

\def\Pr{\mathbb{P}}
\def\le{\leqslant}
\def\ge{\geqslant}
\def\<{\langle}
\def\>{\rangle}
\def\dim{\textup{dim}}

\def\eps{\varepsilon}

\begin{document}
\title{Minimal percolating sets in bootstrap percolation}

\author{Robert Morris}
\address{Murray Edwards College, The University of Cambridge, Cambridge CB3 0DF, England (Work partly done whilst at the University of Memphis)} \email{rdm30@cam.ac.uk}\thanks{The author was supported during this research by a Van Vleet Memorial Doctoral Fellowship.}

\begin{abstract}
In standard bootstrap percolation, a subset $A$ of the grid $[n]^2$ is initially \emph{infected}. A new site is then infected if at least two of its neighbours are infected, and an infected site stays infected forever. The set $A$ is said to \emph{percolate} if eventually the entire grid is infected. A percolating set is said to be minimal if none of its subsets percolate. Answering a question of Bollob\'as, we show that there exists a minimal percolating set of size $4n^2/33 + o(n^2)$, but there does not exist one larger than $(n + 2)^2/6$.
\end{abstract}

\maketitle

\section{Introduction}\label{orderintro}

Consider the following deterministic process on a (finite, connected) graph $G$. Given an initial set of `infected' sites, $A \subset V(G)$, a vertex becomes infected if at least $r \in \N$ of its neighbours are already infected, and infected sites remain infected forever. This process is known as \emph{$r$-neighbour bootstrap percolation} on $G$. If eventually the entire vertex set becomes infected, we say that the set $A$ \emph{percolates} on $G$. For a given graph $G$, we would like to know which sets percolate.

The bootstrap process was introduced in 1979 by Chalupa, Leith and Reich~\cite{CLR}. It is an example of a cellular automaton, and is related to interacting systems of particles; for example, it has been used as a tool in the study of the Ising model at zero-temperature (see~\cite{FSS} and~\cite{me}). For more on the various physical motivations and applications of bootstrap percolation, we refer the reader to the survey article of Adler and Lev~\cite{braz}, and the references therein.

Bootstrap percolation has been extensively studied in the case where $G$ is the $d$-dimensional grid, $[n]^d = \{1,\ldots,n\}^d$, with edges induced by the lattice $\Z^d$, and the elements of the set $A$ are chosen independently at random with probability $p = p(n)$. In particular, much effort has gone into answering the following two questions: $a)$ what is the value of the \emph{critical probability},
$$p_c([n]^d,r) \; = \; \inf\big\{ p \,:\, \Pr_p(A \textup{ percolates}) \ge 1/2 \big\},$$ and $b)$ how fast is the transition from $\Pr(A \textup{  percolates}) = o(1)$ to $\Pr(A \textup{ percolates}) = 1 - o(1)$.

Following fundamental work by Aizenman and Lebowitz~\cite{AL} (in the case $r = 2$) and Cerf and Cirillo~\cite{CC} (in the crucial case $d = r = 3$), Cerf and Manzo~\cite{CM} proved the following theorem, which determines $p_c$ up to a constant for all fixed $d$ and $r$ with $2 \le r \le d$:
$$p_c\big( [n]^d,r \big) \; = \; \Theta\left( \frac{1}{\log_{(r-1)} n} \right)^{d-r+1},$$
where $\log_{(r)}$ is an $r$-times iterated logarithm. Note in particular that $p_c([n]^d,r) = o(1)$ as $n \to \infty$ for every $2 \le r \le d$. More recently much more precise results have been obtained by Holroyd~\cite{Hol}, who proved that in fact
$$p_c\big([n]^2,2\big) \; = \; \frac{\pi^2}{18\log n} \: + \: o\left( \frac{1}{\log n} \right),$$
and by Balogh, Bollob\'as, Duminil-Copin and Morris~\cite{d=r=3,gendr}, who have determined $p_c\big( [n]^d,r \big)$ up to a factor $1 + o(1)$ for all fixed $d$ and $r$. The situation is very different if $d,r \to \infty$ as $n \to \infty$, and there are many open questions. However, very precise results have been obtained by Balogh, Bollob\'as and Morris~\cite{Maj,n^d} (see also~\cite{BB}) in the cases $r = 2$ and $r = d$, as long as $d(n) \to \infty$ sufficiently quickly. For results on other graphs, see~\cite{BPP,BPi,Svante},

As well as studying sets $A \subset [n]^d$ chosen at random, it is very natural to study the extremal properties of percolating sets. For example, it is a folklore fact (and a beautiful exercise to prove) that the minimal size of a percolating set in $[n]^2$ (with $r = 2$) is $n$, and, more generally, the minimal size in $[n]^d$ is $\lceil (n-1)d/2 \rceil + 1$. Perhaps surprisingly, these two questions are closely linked: the lower bound in the result of Aizenman and Lebowitz may be deduced fairly easily from the extremal result, and moreover it is a vital tool in~\cite{n^d}, where the authors determine $p_c([n]^d,2)$ for $d \gg \log n$. Even more surprisingly, the extremal problem is open when $r \ge 3$, even, for example, for the hypercube, $G = [2]^d$. For more results on deterministic aspects of bootstrap percolation, see~\cite{BPe}.

In this paper we shall study a slightly different extremal question, due to Bollob\'as~\cite{Bela}. Given a graph $G$ and a threshold $r$, say that a set $A \subset V(G)$ is a \emph{minimal percolating set} (MinPS) if $A$ percolates in $r$-neighbour bootstrap percolation, but no proper subset of $A$ percolates. Clearly a percolating set of minimal size is a minimal percolating set; but is it true that all minimal percolating sets have roughly the same size? It is the purpose of this note, firstly to introduce the concept of minimal percolating sets, and secondly to show that, contrary to the natural conjecture, there exist fairly dense such sets in $[n]^d$.

We shall study the possible sizes of a minimal percolating set on the $m \times n$ grid, $G(m,n) \subset \Z^2$, with $r = 2$. Let us define
$$E(m,n) = \max\big\{ |A| \,:\, A \subset [m] \times [n]\textup{ is a MinPS of }G(m,n)\big\},$$
and write $E(n) = E(n,n)$. Thus our problem is to determine $E(m,n)$ for every $m,n \in \N$.

It is not hard to construct a minimal percolating set with about $2(m + n)/3$ elements. For example (assuming for simplicity that $m,n \equiv 0 \pmod 3$), take
$$A \; = \; \big\{ (k,1) \,:\, k \equiv 0,2 \textup{ (mod 3)}\} \, \cup \, \{(1,\ell) : \ell \equiv 0,2 \textup{ (mod 3)}\}.$$
It is easy to see that $A$ percolates, and that if $x \in A$, then $A \setminus \{x\}$ does not percolate. For example, if $x = (3,1)$ then the $3^{rd}$ and $4^{th}$ columns of $V = [m] \times [n]$ are empty. However, it is non-trivial to find a MinPS with more than $2(m + n)/3$ elements, and one is easily tempted to suspect that in fact $E(m,n) = \lfloor 2(m + n)/3 \rfloor$. (The interested reader is encouraged to stop at this point and try to construct a minimal percolating set with more than this many elements.) As it turns out, however, the correct answer is rather a long way from this. In fact, even though a randomly chosen set of density $o(1)$ will percolate with high probability, there exist fairly dense minimal percolating sets in $G(m,n)$. The following theorem is the main result of this paper.

\begin{thm}\label{mbyn}
For every $2 \le m,n \in \N$, we have
$$\frac{4mn}{33} \,-\, O\big( m^{3/2} + n\sqrt{m} \big) \; \le \; E(m,n) \; \le \; \frac{(m + 2)(n + 2)}{6}.$$
In particular, $$\frac{4n^2}{33} \,+\, o(n^2) \; \le \; E(n) \; \le \; \frac{(n + 2)^2}{6}.$$
\end{thm}

We remark that Lemma~\ref{justup} (below) gives an explicit lower bound on $E(m,n)$ when $mn$ is small. We suspect that the constant $4/33$ in the lower bound is optimal.

Although we cannot determine $E(m,n)$ asymptotically, we shall at least prove the following theorem, which implies that $E(n) = cn^2 + o(n^2)$, for some constant $c \in [4/33,1/6]$.

\begin{thm}\label{cn^2}
$\ds\lim_{n \to \infty} \frac{E(n)}{n^2}$ exists.
\end{thm}

The rest of the paper is organised as follows. In Section~\ref{basic} we define \emph{corner-avoiding} minimal percolating sets, which will be instrumental in the proofs of Theorems~\ref{mbyn} and \ref{cn^2}, and prove various facts about them, and in Section~\ref{lower} we deduce the lower bound in Theorem~\ref{mbyn}. In Section~\ref{upper} we prove the upper bound in Theorem~\ref{mbyn}, in Section~\ref{limit} we prove Theorem~\ref{cn^2}, and in Section~\ref{ntod} we show how our construction extends to the graph $[n]^d$, and mention some open questions.

\section{Corner-avoiding sets}\label{basic}

Let $m,n \in \N$ and $V = [m] \times [n]$. Given a set $X \subset V$, write $\<X\>$ for the set of points which are eventually infected if the initial set is $X$. If $Y \subset \<X\>$ then we shall say that $X$ \emph{spans} $Y$, and if moreover $Y \subset \<X \cap Y\>$, then we say that $X$ \emph{internally spans} $Y$.

A \emph{rectangle} is a set $$[(a,b),(c,d)] \, := \, \{(x,y) : \, a \le x \le c, \, b \le y \le d\},$$ where $a,b,c,d \in \N$. For any rectangle $R = [(a,b),(c,d)]$, define $$\dim(R) \, := \, (w(R),h(R)) \, := \, (c - a + 1, d - b + 1).$$
Observe that in $G(m,n)$, $\<X\>$ is always a union of rectangles.

The \emph{top-left corner} of $V$ is the rectangle $J_L = [(1,n-1),(2,n)]$ and the \emph{bottom-right corner} of $V$ is the rectangle $J_R = [(m-1,1),(m,2)]$.

\begin{defn}
Call a minimal percolating set $A \subset V$ \emph{corner-avoiding} if whenever $v \in A$, we have $$\<A \setminus \{v\}\> \cap (J_L \cup J_R) = \emptyset,$$ i.e., if the initially infected sites are a (proper) subset of $A$, then the top-left and bottom-right corners remain uninfected.
\end{defn}

Let
$$E_c(m,n) = \max\big\{ |A| : A \subset [m] \times [n]\textup{ is a corner-avoiding MinPS of }G(m,n) \big\}$$
if such sets exist, and let $E_c(m,n) = 0$ otherwise. As before, write $E_c(n) = E_c(n,n)$. Note that the inequality $E_c(m,n) \le E(m,n)$ follows immediately from the definitions.

We start by showing that corner-avoiding minimal percolating sets exist in $G(m,n)$ for certain values of $m$ and $n$.

\vspace{0.1in}
\[ \unit = 1cm \hskip -8\unit
\medline \dl{0}{0}{8}{0} \dl{8}{0}{8}{5} \dl{8}{5}{0}{5} \dl{0}{5}{0}{0}
\thnline \dl{0}{1}{8}{1} \dl{0}{2}{8}{2} \dl{0}{3}{8}{3} \dl{0}{4}{8}{4}
\dl{1}{0}{1}{5} \dl{2}{0}{2}{5} \dl{3}{0}{3}{5} \dl{4}{0}{4}{5} \dl{5}{0}{5}{5} \dl{6}{0}{6}{5} \dl{7}{0}{7}{5}
\varpt{10000} \pt{0.5}{0.5} \pt{0.5}{2.5} \pt{1.5}{2.5} \pt{3.5}{0.5}
\pt{6.5}{2.5} \pt{4.5}{4.5} \pt{7.5}{4.5} \pt{7.5}{2.5}
\point{0.8}{-1}{Figure 1: A corner-avoiding MinPS}
\]
\vspace{0.1in}

Our construction uses the following simple structures. Given a set $A \subset [m] \times [n]$, and integers $k,\ell \in \N$, define
$$A + (k,\ell) \; := \; \{(i,j) \in \N^2 : (i - k,j - \ell) \in A\}.$$
Now, let $P$ be the pair of points $\{(1,1), (1,3)\}$, and for each $k \in \N$ let
$$L(k) \; := \; \bigcup_{i = 0}^{k - 1} \big( P + (0,3i) \big).$$
Furthermore, for each $a,b \in \N$ let $$L(k;a,b) \; := \; L(k) + (a-1,b-1).$$
Observe that $\<L(k;a,b)\> = [(a,b),(a,b + 3k - 1)]$, and that $L(k;a,b)$ is a minimal spanning set for $\<L(k;a,b)\>$.

\begin{lemma}\label{small}
Let $k \in \N$. Then $$E_c(8,3k+2) \ge 4k + 4.$$
\end{lemma}

\begin{proof}
Let
$$A \; = \; L(k) \cup \big\{(2,3k), (4,1), (5,3k+2), (7,3) \big\} \cup \big( L(k) + (7,2) \big)$$
(see Figure 1). Then $A$ is a corner-avoiding minimal percolating set in $[8] \times [3k+2]$, and $|A| = 4k + 4$.
\end{proof}

\begin{rmk}
The bound of Lemma~\ref{small} is connected to the constant $4/33$ in Theorem~\ref{mbyn} in the following way: given a result of the form $E_c(x,yk) \ge zk$, we shall deduce a lower bound of the form
$$E(n) \; \ge \; \ds\frac{z n^2}{(x+3)y}.$$
The $(x + 3)$ term comes from the fact that in Lemma~\ref{2become1}, below, we need to use three extra columns to `connect' two corner-avoiding minimal percolating sets.
\end{rmk}

The next lemma explains our interest in corner-avoiding minimal percolating sets.

\begin{lemma}\label{2become1}
Let $m,m',n,n' \in \N$, and suppose $E_c(m,n) > 0$, $E_c(m',n') > 0$ and $n' \ge n$. Then $$E_c(m+m'+3,n'+2) \: \ge \: E_c(m,n) + E_c(m',n') + 2.$$
\end{lemma}

\begin{proof}
Let $B \subset [m] \times [n]$ and $C \subset [m'] \times [n']$ be corner-avoiding MinPS, with $|B| = E_c(m,n)$ and $|C| = E_c(m',n')$. Note that $B$ and $C$ exist by assumption. Now, let
$$C' \: = \: C + (m+3,2) \: \subset \: [m+m'+3] \times [n'+2],$$ and let
$$A \: = \: B \, \cup \, \{(m+1,1),(m+3,n'+2)\} \, \cup \, C'.$$

\vspace{0.1in}
\[ \unit = 0.25cm \hskip -35\unit
\medline \dl{0}{0}{35}{0} \dl{35}{0}{35}{24} \dl{35}{24}{0}{24} \dl{0}{24}{0}{0}
\medline \dl{11}{0}{11}{18} \dl{11}{18}{0}{18}
\dl{17}{4}{17}{24} \dl{17}{4}{35}{4}
\thnline \dl{7}{0}{7}{4} \dl{7}{4}{11}{4}
\dl{17}{20}{21}{20} \dl{21}{20}{21}{24}
\varpt{8000}
\pt{12.2}{1.2} \pt{15.8}{22.8}
\point{4}{10}{$\<B\>$} \point{24.5}{14}{$\<C'\>$}
\point{9}{-4}{Figure 2: The set $A$}
\]
\vspace{0.1in}

Then $A$ is a corner-avoiding minimal percolating set in $[m+m'+3] \times [n'+2]$ (see Figure 2), and $|A| = E(m,n) + E(m',n') + 2$.
\end{proof}

It is easy to deduce a quadratic lower bound on $E(n)$ from Lemmas~\ref{small} and~\ref{2become1}. However, we shall work harder to obtain what we suspect is an asymptotically sharp lower bound. We begin with a simple application of Lemma~\ref{2become1}.

\begin{lemma}\label{manybecome1}
Let $k,m,n \in \N$. Then
$$E_c(km+3(k-1),n+2(k-1)) \: \ge \: k E_c(m,n).$$
\end{lemma}

\begin{proof}
The proof is by induction on $k$. The result is trivial if $E_c(m,n) = 0$, so assume $E_c(m,n) > 0$. When $k = 1$ we have equality, so suppose $k \ge 2$ and assume the result holds for $k - 1$. Let $m' = (k-1)m + 3(k-2)$ and $n' = n + 2(k-2) \ge n$, so
$$E_c(m',n') \; \ge \; (k-1)E_c(m,n) \; > \; 0$$
by the induction hypothesis. Thus we may apply Lemma~\ref{2become1} to $m$, $n$, $m'$ and $n'$, which gives
\begin{eqnarray*}
E_c(m+m'+3,n'+2) & \ge & E_c\big( m',n' \big) \: + \: E_c\big( m,n \big)\\
& \ge & (k - 1)E_c(m,n) \: + \: E_c\big( m,n \big) \: = \: kE_c\big( m,n \big)
\end{eqnarray*}
as required.
\end{proof}

We shall need one more immediate application of Lemma~\ref{2become1}.

\begin{lemma}\label{recur}
Let $m,n,t \in \N$. Then
$$E_c(2^t(m + 3) - 3,n + 2t) \; \ge \; 2^tE_c(m,n).$$
\end{lemma}

\begin{proof}
The result is immediate if $E(m,n) = 0$, so assume not. Let $g(x) = 2x + 3$ and note that
$$g^t(x) \: = \: 2^t(x + 3) - 3$$ for every $t \in \N$.

We apply Lemma~\ref{2become1} to $E(m,n)$ $t$ times. To be precise, Lemma~\ref{2become1} with $m = m'$ and $n' = n$ gives $E_c(2m + 3,n+2) \ge 2E_c(m,n)$, and hence $$E_c(g^t(m),n+2t) \ge 2^tE_c(m,n).$$ But $g^t(m) = 2^t(m + 3) - 3$, so the result follows.
\end{proof}

\section{A large minimal set}\label{lower}

We now use the results of the previous section to construct a corner-avoiding minimal percolating set in $G(m,n)$ of size $(4/33 + o(1))mn$. The construction will have three stages. First, we use Lemma~\ref{small} to construct a small corner-avoiding minimal percolating set. Then, using Lemma~\ref{manybecome1}, we put about $\sqrt{m}$ of these together to form a long thin minimal percolating set with the right density. Finally we shall use Lemma~\ref{recur} to obtain the desired subset of $G(m,n)$.

We begin with a simple lemma, which we shall need in order to deduce bounds on $E(m,n)$ from those on $E_c(m',n')$. It says that $E(m,n)$ is increasing in both $m$ and $n$.

\begin{lemma}\label{extend}
If $k \le m$ and $\ell \le n$, then $E(k,\ell) \le E(m,n)$.
\end{lemma}

\begin{proof}
By symmetry, it is enough to prove the lemma in the case that $n = \ell$ and $m = k + 1$. So let $A \subset [m-1] \times [n]$ be a MinPS in $G(m-1,n)$, and observe that $(m-1,a) \in A$ for some $a \in [n]$, since $A$ percolates. We claim that one of the sets $B = A \cup \{(m,a)\}$ and $C = A \cup \{(m,a)\} \setminus \{(m-1,a)\}$ is a MinPS for $G(m,n)$.

First suppose that $C \setminus \{u\}$ percolates in $G(m,n)$ for some $u \in C$. Then $A \setminus \{u\}$ must percolate in $G(m-1,n)$, and $u \neq (m,a)$, since $(m,a)$ is the only element of $C$ in column $m$. This contradicts the minimality of $A$.

Note that $B$ percolates in $G(m,n)$, so we may assume that $C$ does not percolate, but $B \setminus \{v\}$ does percolate for some $v \in B$. But $v \notin \{(m-1,a),(m,a)\}$, since $C = B \setminus \{(m-1,a)\}$ doesn't percolate, and $(m,a)$ is the only the only element of $B$ in column $m$. Hence $A \setminus \{v\}$ percolates in $G(m-1,n)$, which contradicts the minimality of $A$. This contradiction completes the proof.
\end{proof}

We can now prove a good bound on $E(m,n)$ in the case that one of $m$ and $n$ is small, say, $m = o(n)$. In the proof of Theorem~\ref{mbyn}, below, we shall apply the first part of Lemma~\ref{justup} with $M \sim \sqrt{m}$ and $N \sim n$.

\begin{lemma}\label{justup}
For every $M,N \in \N$,
$$E_c(11M - 3,3N + 2M) \; \ge \; 4M(N + 1),$$
and hence for every $m,n \in \N$,
$$E(m,n) \; \ge \; 4\left\lfloor \frac{m + 3}{11} \right\rfloor \left\lfloor \frac{n - 2\left\lfloor \frac{m + 3}{11} \right\rfloor + 3}{3} \right\rfloor \; \ge \; \frac{4}{33}\left(mn - \frac{2m^2}{11} - 7n \right).$$
\end{lemma}

\begin{proof}
The first part follows immediately from Lemmas~\ref{small} and~\ref{manybecome1}. Indeed, applying Lemma~\ref{manybecome1} with $m = 8$, $n = 3N + 2$ and $k = M$, we obtain
$$E_c(11M - 3,3N + 2M) \; \ge \; ME_c(8,3N+2) \; \ge \; 4M(N + 1),$$ by Lemma~\ref{small}.

For the second part, let $M = \left\lfloor \frac{m + 3}{11} \right\rfloor$ and $N = \left\lfloor \frac{n - 2M}{3} \right\rfloor$. The result is trivial if $M(N + 1) \le 0$, and if $N = 0$ and $M \ge 1$ then it follows because
$$E(m,n) \; \ge \; \left\lfloor \frac{2(m + n)}{3} \right\rfloor \; \ge \; \frac{4(m+3)}{11},$$ since $m \ge 8$. So assume that $M \ge 1$ and $N \ge 1$, and note that $m \ge 11M - 3$ and $n \ge 3N + 2N$. Thus, by Lemma~\ref{extend},
$$E(m,n) \; \ge \; E_c(11M - 3,3N + 2M),$$
and the result follows by the first part. The final inequality is trivial.
\end{proof}

We are now ready to prove the lower bound in Theorem~\ref{mbyn}.

\begin{proof}[Proof of the lower bound in Theorem~\ref{mbyn}]
We shall prove that
$$E(m,n) \; \ge \; \ds\frac{4mn}{33} \,-\, O\big( m^{3/2} + n\sqrt{m} \big).$$
Assume that $mn$ is sufficiently large, and that $n \ge \sqrt{m}$, since otherwise the result is trivial.

We shall choose positive integers $M$, $N$ and $t$ such that $m \ge 2^t(m'+3)-3$ and $n \ge n' + 2t$, where $m' = 11M - 3$ and $n' = 3N + 2M$. Observe that for such integers, we have
\begin{eqnarray*}
E(m,n) & \ge & E_c(2^t(m' + 3)-3,n'+2t) \; \ge \; 2^tE_c(m',n') \; \ge \; 2^{t+2}M(N+1),
\end{eqnarray*}
by Lemmas~\ref{recur}, \ref{extend} and \ref{justup}.

Indeed, let $t = \left\lceil \ds\frac{\log_2 m }{2} \right\rceil$, $M = \left\lfloor \ds\frac{1}{11} \left( \ds\frac{m+3}{2^t} \right) \right\rfloor$ and $N = \left\lfloor \ds\frac{n - 2t - 2M}{3} \right\rfloor$. Note that $M,N,t \ge 1$, since $n \ge \sqrt{m} \gg 1$, and that $m' = 11M - 3$ and $n' = 3N + 2M$ satisfy the required inequalities. Note also that $2^t \sim \sqrt{m}$, so $M \sim \sqrt{m}$ and $N \sim n$.

Hence,
\begin{eqnarray*}
E(m,n) & \ge & 2^{t+2}M(N+1) \; \ge \; 2^{t+2} \left( \ds\frac{1}{11} \left( \ds\frac{m+3}{2^t} \right) - 1 \right) \left( \ds\frac{n - 2t - 2M}{3} \right)\\
& \ge & \ds\frac{4mn}{33} \,-\, 2^{t+2} n \,-\, (m+3)(M+t) \; = \; \ds\frac{4mn}{33} \,-\, O\big( m^{3/2} + n\sqrt{m} \big),
\end{eqnarray*}
as required.
\end{proof}

\section{An upper bound}\label{upper}

We shall prove the upper bound in Theorem~\ref{mbyn} by induction, using the partial order on vertex sets given by containment. We begin by proving the base cases.

\begin{thm}
Let $n \in \N$. Then
\begin{enumerate}
\item[$(a)$] $E(m,1) = \left\lfloor \frac{2(m + 1)}{3} \right\rfloor$\\[-1ex]

\item[$(b)$] $E(m,2) = \left\lfloor \frac{2(m + 2)}{3} \right\rfloor$\\[-1ex]

\item[$(c)$] $E(m,3) = \left\lfloor \frac{2(m + 3)}{3} \right\rfloor$
\end{enumerate}
\end{thm}

\begin{proof}
The lower bounds are easy, so we shall only prove the upper bounds. In each case, let $A$ be a minimal percolating set. To prove part $(a)$, simply note that $A \subset [m] \times [1]$ can contain at most two out of three consecutive points.

For part $(b)$, observe that if $A \subset [m] \times [2]$ percolates, there must exist $s,t \in [m]$ such that $(s,1),(t,2) \in A$ and $|s - t| \le 1$. Indeed, if no such $s$ and $t$ exist, then $\<\{(k,1) \in A\}\>$ and $\<\{(k,2) \in A\}\>$ are at distance at least 3. There are thus two cases. If $s = t$ then $(i,j) \notin A$ for $i \in \{s-1,s+1\}$, $j \in \{1,2\}$, and $A$ can contain at most two points from any (other) three consecutive columns (else we could remove the middle point). Therefore
$$|A| \le \left\lfloor \frac{2(s - 1)}{3} \right\rfloor + 2 + \left\lfloor \frac{2(m - s)}{3} \right\rfloor \le \left\lfloor \frac{2m + 4}{3} \right\rfloor.$$
If, on the other hand, $s = t + 1$ say, then $A$ contains at most two points from the set $\{(i,j) : i - s \in \{-3,-2,1,2\}, j \in \{1,2\}\}$, and at most two points from any three consecutive columns outside this set. Thus
$$|A| \le \left\lfloor \frac{2(s - 3)}{3} \right\rfloor + 4 + \left\lfloor \frac{2(m - s - 1)}{3} \right\rfloor \le \left\lfloor \frac{2m + 4}{3} \right\rfloor.$$
The reader can easily check that when $s \le 3$ or $s \ge m - 1$, the calculation is exactly the same.

Part $(c)$ requires a little more work, and will be proved by induction on $m$. Observe that the result follows by parts $(a)$ and $(b)$ if $m \le 2$, and that $E(3,3) = E(4,3) = 4$. So let $m \ge 5$, and assume that the result holds for all smaller $m$.

Suppose first that there exists an internally spanned rectangle $R$, with $\dim(R) = (k,3)$, which does not contain either the $(m-1)^{st}$ or the $m^{th}$ column of $V = [m] \times [3]$. Then either $[m-3] \times [3]$ or $[m-2] \times [3]$ must be internally spanned. In the former case, we have
$$|A| \le E(m-3,3) + 2 \le \left\lfloor \frac{2m}{3} \right\rfloor + 2 = \left\lfloor \frac{2(m + 3)}{3} \right\rfloor,$$ while in the latter case we have
$$|A| \le E(m-2,3) + 1 \le \left\lfloor \frac{2(m + 1)}{3} \right\rfloor + 1 \le \left\lfloor \frac{2(m + 3)}{3} \right\rfloor.$$

So assume that no such rectangle $R$ exists (and similarly for the $1^{st}$ and $2^{nd}$ columns of $V$), and observe that there must therefore exist some internally spanned rectangle $T$ with $\dim(T) = (1,2)$ or $(2,2)$. Indeed, if no such rectangle exists then the sets $\<\{(k,1) \in A\}\>$, $\<\{(k,2) \in A\}\>$ and $\<\{(k,3) \in A\}\>$ are (pairwise) at distance at least 3, as in the proof of part $(b)$. Without loss of generality, we may assume (since $m \ge 5$) that $T$ does not intersect either the $(m-1)^{st}$ or the $m^{th}$ column of $V$.

Now, by allowing $T$ to grow one block at a time, we find that either $[m-3] \times [2]$ is internally spanned, or $[m-2] \times [2]$ is internally spanned, or there exists an internally spanned rectangle $T'$, with $\dim(T') = (\ell,2)$ for some $\ell \in [m-4]$, such that $d(A \setminus T', T') \ge 3$. If $[m-2] \times [2]$ is internally spanned, then
$$|A| \le E(m-2,2) + 2 \le \left\lfloor \frac{2m}{3} \right\rfloor + 2 = \left\lfloor \frac{2(m + 3)}{3} \right\rfloor.$$ Also, if $[m-3] \times [2]$ is internally spanned but $[m-2] \times [2]$ is not, then
$$|A| \le E(m-3,2) + 2 \le \left\lfloor \frac{2(m - 1)}{3} \right\rfloor + 2 \le \left\lfloor \frac{2(m + 3)}{3} \right\rfloor,$$ since if $|A \cap [(m-1,1),(m,3)]| \ge 3$, then $[(m-1,1),(m,3)]$ is internally spanned, which contradicts our earlier assumption.

So, without loss of generality, $T' = [\ell] \times [2]$ is internally spanned, and $d(A \setminus T', T') \ge 3$, for some $\ell \in [m-4]$. But then the rectangle $[(\ell+2,2),(m,3)]$ must be internally spanned, since $A$ percolates and there is no internally spanned $k \times 3$ rectangle $R$ in $V$. Thus
\begin{eqnarray*}
|A| & \le & E(\ell,2) + E(m-\ell-1,2) \; \le \; \left\lfloor \frac{2(\ell + 2)}{3} \right\rfloor + \left\lfloor \frac{2(m - \ell + 1)}{3} \right\rfloor \; \le \; \left\lfloor \frac{2(m + 3)}{3} \right\rfloor,
\end{eqnarray*} and so we are done.
\end{proof}

The following corollary is immediate.

\begin{cor}\label{23}
Let $m \in \{2,3\}$, $n \in \N$, then $E(m,n) \le \ds\frac{(m + 2)(n + 2)}{6}$.
\end{cor}

Let $<_R$ be the following partial order on rectangles in $[m] \times [n]$. First, given $a,c \in [m]$ and $b,d \in [n]$, let $(a,b) <_R (c,d)$ if $\min \{m-a,n-b\} \: > \: \min \{m-c,n-d\}$, or $\min \{m-a,n-b\} \: = \: \min \{m-c,n-d\}$ and $\max\{m-a,n-b\} \: > \: \max\{m-c,n-d\}$. Now, given rectangles $S$ and $T$, let $S <_R T$ if and only if $\dim(S) <_R \dim(T)$.

\begin{obs}\label{simple}
If $(p,q) \le_R (k,\ell)$, then $k\ell + p(n - \ell) + q(m - k) \le mn$.
\end{obs}

\begin{proof}
Note that
\begin{eqnarray*}
k\ell + p(n - \ell) + q(m - k) & = & mn + (m - k)(n - \ell) - (m - p)(n - \ell) - (m - k)(n - q).
\end{eqnarray*}
Now, if $p \le k$ then $(m - k)(n - \ell) \le (m - p)(n - \ell)$, while if $p > k$, then $q < \ell$, and so $(m - k)(n - \ell) \le (m - k)(n - q)$. In either case, the result follows.
\end{proof}

We are now ready to prove the upper bound in Theorem~\ref{mbyn}.

\begin{proof}[Proof of the upper bound in Theorem~\ref{mbyn}]
If $2 \le \min\{m,n\} \le 3$ then the result follows by Corollary~\ref{23}, and note that the result also holds if $m = n = 1$ (though it is in general false when $\min\{m,n\} = 1$).

So let $m,n \in \N$, with $m,n \ge 4$, let $A$ be a minimal percolating set in $V = [m] \times [n]$, and assume that if $[p] \times [q] \subsetneq V$ and $p,q \ge 2$, then $E(p,q) \le (p + 2)(q + 2)/6$. We shall show that $|A| \le (m + 2)(n + 2)/6$. In order to aid the reader's understanding, we shall let $a = 1/6$, $b = 1/3$ and $c = 2/3$, so that $(m + 2)(n + 2)/6 = amn + b(m + n) + c$.

Let $S$ be a maximal (in the order $<_R$) internally spanned rectangle in $V$, other than $V$ itself, and let $\dim(S) = (k,\ell)$. We shall distinguish several cases.\\

\noindent Case 1: Either $k = m$ or $\ell = n$.\\

Suppose that $k = m$. Since $A$ percolates, there cannot be two consecutive empty rows, so since $S$ is maximal, we must have $\ell = n - 2$ or $n - 1$, which means that $|A \setminus S| = 1$ (since $A$ is minimal). Hence, by the induction hypothesis,
\begin{eqnarray*}
|A| & \le & E(m,n-1) + 1 \; \le \; am(n - 1) + b(m + n - 1) + c + 1\\[+1ex]
& = & amn + b(m + n) + c - (am + b - 1) \; \le \; amn + b(m + n) + c,
\end{eqnarray*}
since $m \ge 4$, so $am + b \ge 1$. The proof if $\ell = n$ is identical.\\

Assume from now on that $m - k, n - \ell \ge 1$, and let $B = A \setminus S$. Since $S$ is maximal and $A$ is minimal, it follows that either $|B| = 1$, or $d(S,B) \ge 3$. Now let $T = \<B\>$, and note that $T$ is a single rectangle, since some internally spanned rectangle $T_1 \subset T$ must have $d(S,T_1) \le 2$, and thus $\<S \cup T_1\> = V$ by the maximality of $S$. But $A$ is minimal, so we must have $B \setminus T_1 = \emptyset$, and hence $T = T_1$.

Now, since $\< S \cup T \> = V$ and $S$ is maximal (so $\dim(S) \not<_R \dim(T)$), it follows that $S$ contains some corner of $V$. Thus, without loss of generality, we can assume that $S = [k] \times [\ell]$. Say that $S$ and $T$ \emph{overlap rows} if they both contain an element of some row of $V$, and say that they \emph{overlap columns} if they both contain an element of some column.\\

\noindent Case 2: $S$ and $T$ neither overlap rows, nor overlap columns.\\

This means that $T \subset [k+1,m] \times [\ell+1,n]$, and so
\begin{eqnarray*}
|A| & \le & E(k,\ell) + E(m-k,n-\ell)\\ [+1ex]
& \le & a(k\ell + (m - k)(n - \ell)) + b(k + \ell + (m - k) + (n - \ell)) + 2c\\[+1ex]
& = & amn + b(m + n) + c - (a(k(n - \ell) + \ell(m - k)) - c)\\[+1ex]
& \le & amn + b(m + n) + c
\end{eqnarray*}
since $k, m - k \ge 1$, so $k(n - \ell) + \ell(m - k) \ge n \ge 4$, and $4a - c = 0$.\\

So assume, without loss of generality, that $S$ and $T$ overlap columns. Thus $|B| \ge 2$, so $d(S,B) \ge 3$ and $n - \ell \ge 3$. Furthermore, both of the dimensions of $T$ must be at least two, and $B$ is a MinPS for $T$. Thus, (for exactly the same reasons that $S$ and $T$ must exist), there exist disjointly internally spanned rectangles $P$ and $Q$, such that $P,Q \neq T$, but $\<P \cup Q\> = T$ (see Figure $3(a)$). Choose $P$ and $Q$ so that $\min\{|P|,|Q|\}$ is minimal subject to these conditions. Moreover, given $\min\{|P|,|Q|\}$, choose $P$ and $Q$ to each have both dimensions at least two if possible.

Observe first that $d(S,P) \ge 3$ and $d(S,Q) \ge 3$. Indeed, if $d(S,P) \le 2$ then $\<S \cup P\> = V$ (since $S$ is maximal), and so $B \subset P$ (since $A$ is minimal), so $P = T$. But we chose $P \neq T$, so this is a contradiction. Let $\dim(P) = (p,s)$ and $\dim(Q) = (t,q)$.

\vspace{0.1in}
\[ \unit = 0.6cm \hskip -19\unit
\medline \dl{0}{0}{8}{0} \dl{8}{0}{8}{8} \dl{8}{8}{0}{8} \dl{0}{8}{0}{0}
\dl{5}{0}{5}{5} \dl{5}{5}{0}{5} \dl{2.5}{8}{2.5}{4} \dl{2.5}{4}{8}{4}
\dl{2.5}{6}{7}{6} \dl{7}{6}{7}{8} \dl{8}{7.5}{6}{7.5} \dl{6}{7.5}{6}{4}
\point{2}{2}{\small $S$} 
\point{4}{6.8}{\small $P$} \point{7}{5}{\small $Q$}
\point{2}{6.8}{\tiny $s$} \point{8.2}{5.7}{\tiny $q$} \point{4.5}{8.2}{\tiny $p$} \point{7}{3.5}{\tiny $t$}
\point{2.2}{-0.5}{\tiny $k$} \point{-0.5}{2.3}{\tiny $\ell$}
\dl{12}{4}{17}{4} \dl{17}{4}{17}{1} \point{14}{2}{\small $S$}
\dl{16}{5}{16}{8} \dl{12}{8}{19}{8} \dl{19}{8}{19}{1} \dl{19}{5}{16}{5} \point{15}{6.5}{\small $T$}
\varpt{10000} \pt{16.5}{6.5} \pt{18.5}{5.5} \pt{18.5}{7.5}
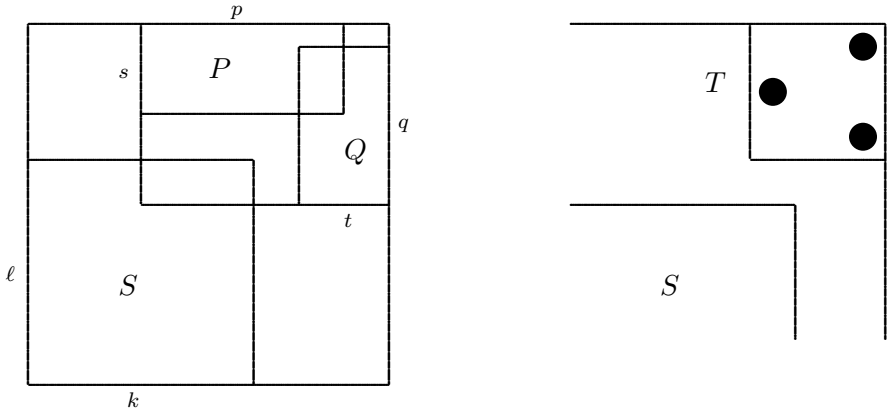
\point{-0.8}{-2}{\small Figure 3: $(a)$ The rectangles $S$, $P$ and $Q$, $(b)$ a configuration with $|B| = 3$}
\]
\vspace{0.1in}

We claim that either $\min\{p,q,s,t\} \ge 2$, or $|Q| = 1$ and both dimensions of $P$ are at least two (or vice-versa), or $|B| \le 3$. Note than in the first two cases we can apply the induction hypothesis to $P$ and $Q$; the third case is illustrated in Figure $3(b)$.

Suppose first that $\min\{q,t\} = 1$ but $|B \cap Q| \ge 3$. Let $u_1$ and $u_2$ be the end-vertices of $Q$, and let $Q_i = \< B \cap Q \setminus u_i\>$ for $i = 1,2$. Since $d(P,Q) \le 2$ and $|B \cap Q| \ge 3$, it follows that $d(P,Q_i) \le 2$ for some $i \in \{1,2\}$. Therefore we could have chosen $P$ and $Q$ with $|Q| = 1$ and both dimensions of $P$ at least two; in particular $P^* = \<P \cup Q_i\>$ and $Q^* = u_i$ would do. (Note that $P^* \neq T$ since $A$ is minimal.) This contradicts our choice of $P$ and $Q$. Similarly, we cannot have $\min\{s,p\} = 1$ and $|B \cap P| \ge 3$.

So suppose that $\min\{q,t\} = 1$ and $\max\{q,t\} \in \{2,3\}$. If $|P| = 1$ then $|B| = 3$ (see Figure $3(b)$). But if $\min\{s,p\} \ge 2$, then $d(P,u) \le 2$ for some $u \in B \cap Q$ (since $d(P,Q) \le 2$), and we could replace the pair $(P,Q)$ by $(P^*,Q^*) = (\< P \cup Q_i\>,u_i)$, as above.

Thus we may assume that $\min\{s,p\} = \min\{q,t\} = 1$, and that $B \cap Q = \{u_1,u_2\}$ and $B \cap P = \{v_1,v_2\}$. But now, again using the fact that $d(P,Q) \le 2$, it is easy to see that either $d(P,u_i) \le 2$ or $d(Q,v_i) \le 2$ for some $i \in \{1,2\}$. Therefore we could have chosen $P$ and $Q$ with $|Q| = 1$, and we have our final contradiction.

We conclude that either $\min\{p,q,s,t\} \ge 2$, or (without loss of generality) $|Q| = 1$ and both dimensions of $P$ are at least two, or $|B| \le 3$, as claimed.\\

\noindent Case 3: $|B| \le 3$.\\

Recall that $m - k \ge 1$ and $n - \ell \ge 3$, so
\begin{eqnarray*}
|A| & \le & E(m-1,n-3) + 3\\[+1ex]
& \le & a(m-1)(n-3) + b(m + n - 4) + c + 3\\[+1ex]
& = & amn + b(m + n) + c - (a(3m + n - 3) + 4b - 3)\\[+1ex]
& \le & amn + b(m + n) + c,
\end{eqnarray*}
since $m,n \ge 4$, so $a(3m + n - 3) + 4b - 3 \ge 13a + 4b - 3 > 0$.\\

So assume from now on that $|B| \ge 4$, and so by the comments above, both dimensions of $P$ are at least two, and either $|Q| = 1$ or both dimensions of $Q$ are least two. In particular, we may apply the induction hypothesis to the rectangles $P$ and $Q$.\\

\noindent Case 4: $|Q| = 1$, and $S$ and $P$ overlap columns.\\

We shall need the following simple inequality:
$$an(m - p) \,+\, b(m - k) \; \ge \; a\ell(k - p) \,+\, 1.$$
To see this, suppose first that $k \ge p$, so $a\ell(k - p) \le an(k - p)$. The inequality is thus implied by $(m-k)(an + b) \ge 1$, which holds because $m - k \ge 1$ and $an + b = 1$.

On the other hand, if $k < p$ then $a\ell(k - p) < 0$. But $m - p \ge 1$ (since $S$ was maximal in the order $<_R$ and $m - k \ge 1$), so the inequality follows from the fact that $an + b = 1$.

Now, since $d(S,P) \ge 3$ and $S$ and $P$ overlap columns, we have $s \le n - \ell - 2$. It follows that
\begin{eqnarray*}
|A| & \le & E(k,\ell) + E(p,n-\ell-2) + 1\\ [+1ex]
& \le & a(k\ell + pn - p\ell - 2p)) + b(k + p + n - 2) + 2c + 1\\[+1ex]
& = & a(pn + \ell(k - p)) + b(k + n) - p(2a - b) - (2b - c) + c + 1\\[+1ex]
& = & a(mn  - n(m - p) + \ell(k - p)) + b(m + n - (m - k)) + c + 1\\[+1ex]
& \le & amn + b(m + n) + c
\end{eqnarray*}
by the inequality above, and since $2a = b$ and $2b = c$.\\

\noindent Case 5: $|Q| = 1$, and $S$ and $P$ do not overlap columns.\\

First note that
$$(k + 1)(n - \ell) \,+\, \ell(m - k) \; \ge \; m + 2k + 3 \; \ge \; 9,$$ since $m \ge 4$, $k,\ell \ge 1$ and $n - \ell \ge 3$. Similarly $k(n - \ell) + (\ell + 1)(m - k) \ge 2m + k \ge 9$. Note also that $9a + b - c - 1 \ge 0$.

Now, since $S$ and $T$ overlap columns, $S$ and $Q$ must overlap columns. But $|Q| = 1$, $d(P,Q) \le 2$ and $d(S,P) \ge 3$, so $S$ and $P$ cannot overlap rows. Also $(k+1,\ell+1) \not\in P$. Thus, by the inequalities above, either
\begin{eqnarray*}
|A| & \le & E(k,\ell) + E(m-k-1,n-\ell) + 1\\ [+1ex]
& \le & a(k\ell + (m - k - 1)(n - \ell)) + b(m + n - 1) + 2c + 1\\[+1ex]
& = & amn + b(m + n) + c - (a((k + 1)(n - \ell) + \ell(m - k)) + b - c - 1)\\[+1ex]
& \le & amn + b(m + n) + c,
\end{eqnarray*} or
\begin{eqnarray*}
|A| & \le & E(k,\ell) + E(m-k,n-\ell-1) + 1\\ [+1ex]
& \le & amn + b(m + n) + c - (a(k(n - \ell) + (\ell + 1)(m - k)) + b - c - 1)\\[+1ex]
& \le & amn + b(m + n) + c,
\end{eqnarray*}
as required.\\

So assume from now on that $\min\{p,q,s,t\} \ge 2$. Since $S$ and $T$ overlap columns, we may assume without loss of generality that $S$ and $P$ overlap columns. Since $d(S,P) \ge 3$ and $s \ge 2$, it follows that $n - \ell \ge 4$.

Recall that $\dim(T) = (w(T),h(T))$, and note that $h(T) \ge q + 1$, since otherwise we could have chosen $P$ and $Q$ with $|P| = 1$. Similarly, $w(T) \ge p + 1$.\\

\noindent Case 6: $S$ and $T$ overlap both rows and columns.\\

Since $d(S,P) \ge 3$ and we assumed that $S$ and $P$ overlap columns, it follows that $s \le n - \ell - 2$ and $S$ and $P$ do not overlap rows. Therefore $S$ and $Q$ must overlap rows, and so $t \le m - k - 2$. Note also that $(p,q) \le_R \dim(T) \le_R (k,\ell)$, so we may apply Observation~\ref{simple}. Hence,
\begin{eqnarray*}
|A| & \le & E(k,\ell) + E(p,n-\ell-2) + E(q,m-k-2)\\ [+1ex]
& \le & a(k\ell + p(n - \ell - 2) + q(m - k - 2)) + b(m + n + p + q - 4) + 3c\\[+1ex]
& \le & amn + b(m + n) - (p + q)(2a - b) - 2(2b - c) + c\\[+1ex]
& = & amn + b(m + n) + c
\end{eqnarray*}
since $2a = b$ and $2b = c$, and by Observation~\ref{simple}.\\

There is one remaining case to consider.\\

\noindent Case 7: $S$ and $T$ overlap columns but not rows.\\

Let
$$M \; := \; (k - p)(n - \ell) + (m - k)(n - q) \; = \; (m - p)(n - \ell) + (m - k)(\ell - q).$$
We shall use the following two facts about $M$:
\begin{itemize}
\item $mn - M \; = \; k\ell + p(n - \ell) + q(m - k)$.\\[-1.5ex]
\item $M \; \ge \; 2q + 4$.
\end{itemize}
The first fact is straightforward. To see the second, first suppose that $k \ge p + 2$, and note that $q \le n - \ell - 1$, since $h(T) \ge q + 1$, and $p \le m - 2$, since $w(T) \ge p + 1$. Then
$$M \; \ge \; (k - p)(n - \ell) + (m - k)(\ell + 1) \; \ge \; 2(n - \ell) + 2 \; \ge \; 2q + 4.$$
But if $k \le p + 1 \le w(T)$, then we must have $\ell \ge h(T) \ge q - 1$, since $S \not<_R T$. If $k \le m - 2$ then we are now done, since
$$M \; \ge \; (m - p)(n - \ell) + (m - k) \; \ge \; 2(n - \ell) + 2 \; \ge \; 2q + 4.$$
But if $k = m - 1$ then $B$ contains no element of the $(\ell + 1)^{st}$ row, since $d(S,B) \ge 3$. It follows that $q \le h(T) - 1 \le n - \ell - 2$, and so
$$M \; \ge \; (m - p)(n - \ell) \; \ge \; 2(n - \ell) \; \ge \; 2q + 4,$$
as required.

Now, recall that $d(S,P) \ge 3$ and that $S$ and $P$ overlap columns, so $s \le n - \ell - 2$. Recall also that $d(S,Q) \ge 3$ and $d(S,T) \le 2$, so $S$ and $Q$ do not overlap columns, and so $t \le m - k$. Hence, using the two facts proved above, we have
\begin{eqnarray*}
|A| & \le & E(k,\ell) + E(p,n-\ell-2) + E(q,m-k)\\ [+1ex]
& \le & a(k\ell + p(n - \ell - 2) + q(m - k)) + b(m + n + p + q - 2) + 3c\\[+1ex]
& = & a(mn - M - 2p) + b(m + n) - b(p + q) - 2(b - c) + c\\[+1ex]
& \le & amn + b(m + n) - (p + q)(2a - b) - 2(2a + b - c) + c\\[+1ex]
& = & amn + b(m + n) + c,
\end{eqnarray*}
since $2a = b$ and $2a + b = c$, and we are done.
\end{proof}

\section{Proof of Theorem~\ref{cn^2}}\label{limit}

In this section we shall prove that the sequence $E(n)/n^2$ converges. The proof uses Lemma~\ref{2become1}, together with the following, probably well-known, result on (almost) super-additive sequences: it is a two-dimensional version of Fekete's Lemma. For completeness we shall sketch the proof.

\begin{lemma}\label{2dFek}
Let $M \in \N$, and suppose $f : \N \times \N \to \N$ satisfies $f(m,n) = f(n,m)$, and
\begin{equation} \label{eq1} f(m + m' + 3, n' + 2) \ge f(m,n) + f(m',n') \end{equation}
for every $M \le m,m',n,n' \in \N$ with $n \le n'$. Then $\ds\frac{f(n,n)}{n^2}$ converges as $n \to \infty$.
\end{lemma}

\begin{proof}
For simplicity, we shall write $f(n) = f(n,n)$, and assume that $M$ is large. Let $\eps > 0$, and suppose $f(k) \ge ck^2$ for some sufficiently large $k \in \N$ and some $c \in [0,1]$ (in particular let $k \gg M^2$).

We shall prove that $f(n) \ge (c - \eps)n^2$ for every sufficiently large $n \in \N$. This follows from the following three claims.

\bigskip

\noindent \ul{Claim 1}: $f(n,n) \ge f(m_1,m_2)$ if $n \ge m_i + M + 3 \ge 2M + 3$ for each $i \in \{1,2\}$.

\begin{proof}[Proof of claim]
\begin{eqnarray*}
f(n,n) \; \ge \; f(m_1,m_2) \, + \, f(n - m_1 - 3, n - 2) \; \ge \; f(m_1,m_2),
\end{eqnarray*}
as required.
\end{proof}

\bigskip

\noindent \ul{Claim 2}: $f(tm + 4t^2) \ge t^2f(m)$ for every $m,t \ge M$.

\begin{proof}[Proof of claim]
Applying inequality \eqref{eq1} $t - 1$ times, we obtain
$$f(m + 2(t-1), tm + 3(t-1)) \; \ge \; t f(m,m),$$
and similarly
$$f(tm + 2t(t-1) + 3(t - 1), tm + 5(t-1)) \; \ge \; t f(m + 2(t - 1), tm + 3(t - 1)).$$
Hence, by Claim 1,
$$f(tm + 4t^2) \; \ge \; f(tm + 2t^2 + t - 3, tm + 5(t-1)) \; \ge \; t^2f(m,m)$$
as required.
\end{proof}

\bigskip

\noindent \ul{Claim 3}: $f\big(2^r(m + 2M) - 2M\big) \ge 4^r f(m)$ for every $r \ge 1$ and $m \ge M$.

\begin{proof}[Proof of claim]
We have, by Claim 1 and inequality \eqref{eq1},
$$f(2m + 2M) \; \ge \; f(2m + 7,2m+5) \; \ge \; 2f(m + 2, 2m + 3) \; \ge \; 4f(m,m).$$
Iterating $r$ times, we obtain the required inequality.
\end{proof}

Now, let
$$n' \; = \; 2^r\big( tk + 4t^2 + 2M \big) \, - \, 2M,$$
where $r$ is chosen so that $k^{3/2} \le \ds\frac{n}{2^r} \le 2k^{3/2}$, and $t$ is chosen so that
$$\frac{n}{2^r} \, - \, 2k \; \le \; tk + 4t^2 + 2M \; \le \; \frac{n}{2^r}.$$
Then $n' \le n - 2M$ and $t = \Theta(\sqrt{k})$, and so, by Claims 1, 2 and 3,
\begin{eqnarray*}
f(n) \; \ge \; f(n') \; \ge \; 4^r f(tk + 4t^2) & \ge & 4^r t^2 f(k) \; \ge \; t^2 f(k) \left( \frac{n}{tk + 4t^2 + 2M + 2k} \right)^2 \\
& \ge & \frac{f(k) n^2}{(k + O(\sqrt{k}))^2} \; \ge \; (c - \eps)n^2
\end{eqnarray*}
if $n$ and $k$ are sufficiently large, as required.
\end{proof}

It follows easily from Lemmas~\ref{2become1} and \ref{2dFek} that $E_c(n)/n^2$ converges as $n \to \infty$. However, in order to show that $E(n)/n^2$ also converges we need the following simple result, which follows from the techniques of Section~\ref{basic}.

\begin{lemma}\label{cavreg}
If $m,n \in \N$, with $n \ge 4$, then $$E_c(m,n) \; \le \; E(m,n) \; \le \; E_c(m+16,n+8) \, - \, \frac{4n}{3}.$$
\end{lemma}

\begin{proof}[Proof of Lemma~\ref{cavreg}]
The first inequality is obvious from the definition; we shall prove the second inequality. Let $m,n \in \N$, with $n \ge 4$, and let $A$ be a minimal percolating set of $G(m,n)$ with $|A| = E(m,n)$.

Recall from Section~\ref{basic} the definition of $L(k)$. Let $m' = m + 16$ and $n' = n + 8$, let $V = [m'] \times [n']$, let $N = \lfloor \frac{n + 2}{3} \rfloor$, and define
$$B = L(N) \cup \{(2,3N), (4,1), (5,n+4), (7,5)\}.$$ Now, let $B'$ be the set obtained by rotating $B$ through $180^\circ$, and placing the top-right corner at the point $(m',n')$, i.e.,
$$B' = \{(x,y) : (m' - x + 1,n' - y + 1) \in B\}.$$
Finally, let
$$C = B \cup (A + (8,4)) \cup B',$$ so $C \subset V$ (see Figure 4).

\[ \unit = 0.32cm \hskip -31\unit
\medline \dl{0}{0}{30}{0} \dl{30}{0}{30}{20} \dl{30}{20}{0}{20} \dl{0}{20}{0}{0}
\varpt{4000} \pt{0.5}{0.5} \pt{0.5}{2}
\point{0}{2.25} {$ \pt{0.5}{0.5} \pt{0.5}{2} $}
\medline \dline{0.5}{4.5}{0.5}{9}{4}
\point{0}{9} {$ \pt{0.5}{0.5} \pt{0.5}{2} $}
\point{0}{11.25} {$  \pt{0.5}{0.5} \pt{0.5}{2} $}
\point{1}{12.75} {$  \pt{0.5}{0.5} $}
\point{3}{0} {$  \pt{0.5}{0.5} $}
\point{0}{-1.95} {$ \point{4}{17.9} {$  \pt{0.5}{0.5} $} \point{5.8}{5} {$  \pt{0.5}{0.5} $}
\medline \dl{7.8}{5}{7.8}{18.9} \dl{7.8}{18.9}{22.2}{18.9} \dl{22.2}{18.9}{22.2}{5}\dl{22.2}{5}{7.8}{5} $}
\point{11.5}{9.5}{\small $\<A + (8,4)\>$}
\point{30}{20} {$ \varpt{4000} \pt{-0.5}{-0.5} \pt{-0.5}{-2}
\point{-0}{-2.25} {$ \pt{-0.5}{-0.5} \pt{-0.5}{-2} $}
\medline \dline{-0.5}{-4.5}{-0.5}{-9}{4}
\point{-0}{-9} {$ \pt{-0.5}{-0.5} \pt{-0.5}{-2} $}
\point{-0}{-11.25} {$  \pt{-0.5}{-0.5} \pt{-0.5}{-2} $}
\point{-1}{-12.75} {$  \pt{-0.5}{-0.5} $}
\point{-3}{-0} {$  \pt{-0.5}{-0.5} $}
\point{-0}{1.95} {$ \point{-4}{-17.9} {$  \pt{-0.5}{-0.5} $} \point{-5.8}{-5} {$  \pt{-0.5}{-0.5} $} $} $}
\point{8}{-3}{Figure 4: The set $C$}
\point{4}{0.5}{\tiny $(4,1)$} \point{1.8}{13.5}{\tiny $(2,3N)$} \point{0.4}{16.8}{\tiny $(5,n+4)$} \point{3.8}{3.7}{\tiny $(7,5)$}
\]
\vspace{0.1in}

It is easy to see that $C$ is a corner-avoiding minimal percolating set of $G(m',n')$. Since $|B| \ge 2N \ge 2n/3$, it follows immediately that
$$E_c(m+16,n+8) \; \ge \; |C|  \; \ge \; |A| + \frac{4n}{3} \; = \; E(m,n) + \frac{4n}{3},$$ as required.
\end{proof}

Finally, we may deduce Theorem~\ref{cn^2}.

\begin{proof}[Proof of Theorem~\ref{cn^2}]
By Lemma~\ref{2become1}, we have
$$E_c(m + m' + 3, n' + 2) \: \ge \: E_c(m,n) + E_c(m',n') + 2$$
for every $m,m',n,n' \in \N$ such that $E_c(m,n) > 0$, $E_c(m',n') > 0$ and $n' \ge n$. Since $E_c(8,3k+2) > 0$ for every $k \in \N$, by Lemma~\ref{small}, it is straightforward to deduce that $E_c(m,n) > 0$ if $m$ and $n$ are sufficiently large. Thus, by Lemmas~\ref{2dFek} and \ref{cavreg},
$$\lim_{n \to \infty} \frac{E(n)}{n^2} \; = \; \lim_{n \to \infty} \frac{E_c(n)}{n^2}$$
exists, as required.
\end{proof}

\section{Further problems}\label{ntod}

We have been studying a special case of a much more general question. Indeed, for each graph $G$, and each $r \in \N$, we may define
\begin{eqnarray*}
E(G,r) & := & \max\{ \, |A| \, : \, A \subset V(G)\textup{ is a minimal percolating set of }G \\
&& \hspace{6cm} \textup{in } r \textup{-neighbour bootstrap percolation} \}.
\end{eqnarray*}

\begin{prob}
Determine $E(G,r)$ for every graph $G$ and $2 \le r \in \N$. In particular, does there exist a bounded degree graph sequence $(G_n)$ such that $E(G_n,r) = o(n)$?
\end{prob}

The following straightforward corollary of Theorem~\ref{mbyn} shows that $[n]^d$ is not such a graph sequence for $r = 2$.

\begin{thm}
There exists a function $C: \N \to \N$ such that $$E([n]^d,2) \; \ge \; C(d)n^d,$$
for every $n,d \in \N$.
\end{thm}

\begin{proof}[Sketch of proof]
Divide $[n]^d$ into $n$ hyperplanes, each of co-dimension 1. Take a (corner-avoiding) construction for $[n]^{d-1}$ in every fourth hyperplane, and put a single point in opposite corners of every first and third (mod 4) hyperplane, along the lines of Lemma~\ref{2become1}. This set is now a corner-avoiding minimal percolating set of $[n]^d$.
\end{proof}

We finish with a question, and a conjecture.

\begin{prob}
What is the behaviour of $\:\ds\frac{E([n]^d,2)}{n^d}\:$ as $d = d(n) \to \infty$?
\end{prob}

\begin{conj}
$\ds\lim_{n \to \infty} \frac{E(n)}{n^2} \, = \, \ds\frac{4}{33}.$
\end{conj}

\end{document}